\documentclass{amsart}

\usepackage[utf8]{inputenc}
\usepackage[T1]{fontenc}
\usepackage{eulervm}
\usepackage{tgpagella}

\usepackage[numbers]{natbib}

\usepackage{amsmath}
\usepackage{amssymb}
\usepackage{amsthm}
\usepackage{amsfonts}
\usepackage{mathrsfs}
\usepackage{xspace}
\usepackage{tikz}
\usepackage{nicefrac}
\usepackage{fixmath}
\usepackage{paralist}
\usepackage{ellipsis}
\usepackage{mathtools}
\usepackage{graphicx,subfigure,epic,eepic}
\usepackage{array}

\usepackage{thmtools}
\usepackage{thm-restate}

\usepackage{cleveref}

\usepackage{MnSymbol}

\theoremstyle{definition}
\newtheorem{definition}{Definition}[section]

\theoremstyle{plain}
\declaretheorem[name=Theorem, numberwithin=section]{theorem}
\declaretheorem[numberlike=theorem,name=Lemma]{lemma}
\declaretheorem[numberlike=theorem,name=Proposition]{proposition}
\declaretheorem[numberlike=theorem,name=Corollary]{corollary}

\theoremstyle{remark}
\declaretheorem[name=Remark,numberwithin=section]{remark}
\declaretheorem[name=Example,numberwithin=section]{example}

\DeclareMathOperator{\tr}{Tr}

\DeclareMathOperator{\id}{Id}

\DeclareMathOperator{\rk}{rk}
\DeclareMathOperator{\Pic}{\underline{Pic}}
\DeclareMathOperator{\pic}{Pic}
\DeclareMathOperator{\Br}{Br}

\DeclareMathOperator{\Sym}{Sym}
\DeclareMathOperator{\Spec}{Spec}
\DeclareMathOperator{\NS}{NS}

\newcommand{\Z}{\mathbb{Z}\xspace}

\newcommand{\cM}{\mathcal{M}\xspace}

\newcommand{\C}{\mathbb{C}\xspace}

\newcommand{\cF}{\mathcal{F}\xspace}

\newcommand{\bP}{\mathbb{P}\xspace}
\newcommand{\cY}{\mathcal{Y}\xspace}
\newcommand{\cA}{\mathcal{A}\xspace}
\newcommand{\cO}{\mathcal{O}\xspace}

\newcommand{\sD}{\mathscr{D}\xspace}

\title{Mordell-Weil Groups of Linear Systems and the Hitchin Fibration}
\author{Matthew Woolf}
\begin{document}

\begin{abstract}
In this paper, we study rational sections of the relative Picard scheme of a linear system on a smooth projective variety. We prove that if the linear system is basepoint-free and the locus of non-integral divisors has codimension at least two, then all rational sections of the relative Picard scheme come from restrictions of line bundles on the variety. As a consequence, we describe the group of sections of the Hitchin fibration for moduli spaces of Higgs bundles on curves.\end{abstract}

\maketitle

\section{Introduction}

Given a family of smooth projective varieties $\cY \to S$, the relative Picard scheme is a family of abelian varieties $\cA \to S$. The generic fiber $\cA_\eta$ of $\cA$ is an abelian variety over a non-closed field $k(\eta)$. The $k(\eta)$-rational points of $A_\eta$ correspond to rational sections of the fibration $\cA \to S$ and form an abelian group, called the \emph{Mordell-Weil group}. One often expects that the Mordell-Weil group is as small as possible. In this paper, we study the Mordell-Weil group when $\cY \to S$ is a linear system of divisors satisfying mild assumptions on a smooth projective variety $X$. We verify the expectation in this case by showing that the Mordell-Weil group is generated by line bundles on $X$. 

The inspiration for this paper is the strong Franchetta conjecture. The canonical bundle gives a natural section of the universal Picard scheme over the moduli space of curves. The strong Franchetta conjecture, proved by Mestrano \citep{mestrano} and Kouvidakis \citep{kouvidakis}, states that this section generates the Mordell-Weil group.

Let $|D|$ be a basepoint-free linear system of divisors in a smooth projective variety $X$ with universal family $\sD \to |D|$. Any line bundle $L$ on $X$ gives a rational section of the relative Picard scheme $\Pic(\sD/|D|)$ (see Remark \ref{singularpicard}). Our main theorem verifies that every element of the Mordell-Weil group is of this form, provided the non-integral divisors in $|D|$ have codimension at least two.

\begin{restatable}{theorem}{bigthmone}\label{mainresult}
Let $X$ be a smooth projective variety over $\C$, $|D|$ a basepoint-free linear system of divisors on $X$ such that the locus of divisors which are not integral has codimension at least two. Then rational sections of the relative Picard scheme all come from restricting line bundles on $X$.
\end{restatable}

This theorem can be thought of as a variant of the Grothendieck-Lefschetz theorem asserting that if $D \subset X$ is a smooth ample divisor of dimension at least three, then the restriction map $\pic(X) \to \pic(D)$ is an isomorphism. We make no assumptions about ampleness or dimension, but we only get a result about the relative Picard group, not the individual Picard groups of each smooth divisor in the linear system. In the case of curves, where the Grothendieck-Lefschetz theorem does not apply, our result is much more useful. For example, by considering curves in K3 surfaces of Picard rank 1 which generate the Picard group, one can show that the degree of any rational section of the universal Picard scheme over the moduli space of curves must have degree divisible by $2g-2$. %Using results from \citep{mestrano}, one may then deduce the original strong Franchetta conjecture.% We can for example immediately deduce that on a surface in $\bP^3$ of Picard number one and degree at least four, the general curve in any linear system on the surface is not hyperelliptic.

The hypothesis on the dimension of the locus of non-integral divisors is necessary. We will provide counterexamples without this hypothesis (see Section \ref{counterexamples}). Moreover, by the Castelnuovo-Kronecker theorem \citep{Castelnuovo}, this hypothesis holds as long as the image of the map given by $|D|$ is not ruled by lines or a projection of the Veronese surface. %Moreover, in the case of very ample linear systems on a surface, the Castelnuovo-Kronecker theorem (see \citep{Castelnuovo}) implies that this hypothesis is satisfied except in the case where the image of the map from the surface to projective space is either ruled by lines or a projection of the Veronese surface.

We apply Theorem \ref{mainresult} to study sections of the Hitchin fibration on the moduli space of Higgs bundles on a curve $\Sigma$. One can interepret this fibration birationally as the relative Picard scheme of a linear system on the total space of the cotangent bundle $T^*\Sigma$. The Hitchin fibration has natural sections coming from line bundles on $T^*\Sigma$. %Another case where we have a family of abelian varieties comes from considering moduli spaces of stable Higgs bundles. These moduli spaces admit a map to affine space, called the Hitchin fibration, with general fiber an Abelian variety. It is possible to produce some natural rational sections of the Hitchin fibration. Another result of this paper is that almost always, these are the only sections.

\begin{restatable}{theorem}{bigthmtwo}\label{hitchinresult}
Let $\Sigma$ be a curve of genus $g$, and $M(r,d)$ the moduli space of rank $r$ Higgs bundles of degree $d$. If $g>2$ and $r>1$ or $g=2$ and $r>2$, then all sections of the Hitchin fibration are given by lines bundles on $T^*\Sigma$.
\end{restatable}

%Our proof strategy will be to first prove the analogue of the weak Franchetta conjecture, which is very easy in this setting, then to restrict our section to pencils, where we can use Tsen's theorem to show that any rational section must be a linear combination of the restriction of a line bundle on $X$ and the base locus of the pencil. We then show that in fact we can choose a line bundle on $X$ which gives rise to the section.

I would like to thank Dawei Chen and Joe Harris for their many helpful conversations on this topic. I would also like to thank Arend Bayer, Izzet Coskun, Clifford Earle, Nicole Mestrano, Brendan Hassett, Steve Kleiman, Eric Riedl, and Laura Schaposnik for their help.

We work over $\mathbb{C}$ for simplicity. %For us, a curve will be a connected reduced projective scheme of dimension one.

\section{Relative Picard Varieties}

In this section, we collect the necessary details about relative Picard varieties and their sections. For further details, we refer the reader to \citep{kleiman} for relative Picard varieties and \citep{milne} for Brauer groups.

Given a scheme $Y$, let $\pic(Y)$ be the group of line bundles on $Y$. Given a morphism $f:Z \to Y$, let $\pic(Z/Y)$ be $\pic(Z)/f^*\pic(Y)$.

Given a smooth projective morphism of varieties $\pi:\sD \to S$ with geometrically connected fibers, \citep[Theorem 4.8]{kleiman} proves the existence of an associated relative Picard scheme, denoted by $\Pic(\sD/S)$ with a natural map $p:\Pic(\sD/S) \to S$. The fiber of $p$ over a geometric point $s \in S$ is the Picard scheme $\Pic(\sD_s)$. The relative Picard scheme is a countable disjoint union of projective varieties indexed by the relative Neron-Severi group, $\NS(\sD/S)$. Given an element $\tau \in \NS(\sD/S)$, we define $\Pic^\tau(\sD/S)$ to be the corresponding component of the relative Picard scheme.

For the rest of this section, $S$ is a smooth connected variety and $\pi:\sD \to S$ is a smooth projective morphism with geometrically connected fibers.

\begin{proposition} \label{relpic}
Suppose $L \in \pic(\sD)$ restricts to the trivial line bundle on each fiber of $\pi$. Then $L \cong \pi^* L'$ for some $L' \in \pic(S)$.
\end{proposition}
\begin{proof}
 By Grauert's theorem \citep[III.12.9]{hartshorne}, $L'= \pi_* L$ is a line bundle. The map \[\pi^* L'=\pi^* \pi_* L \to L\] is a nonzero map of invertible sheaves which is an isomorphism on the fibers of $\pi$, again by Grauert's theorem. 
Therefore, the map is an isomorphism of sheaves.
\end{proof}

%This has the following consequence.

The following proposition is \citep[Corollary 1.5]{mestranoramanan} if the fibers of $\pi$ are curves. The same proof works in general.

\begin{proposition} \label{section_extend}
Let $\tau \in \mathrm{NS}(\sD/S)$. Let $\sigma:S \to \Pic^\tau(\sD/S)$ be a rational section of the natural map $\Pic^\tau(\sD/S) \to S$. Then $\sigma$ extends to a regular section.
\end{proposition}

This proposition allows us to ignore the distinction between regular sections and rational sections of the relative Picard scheme.

\begin{remark}\label{singularpicard}
If $\pi:\sD/S$ is only generically smooth, and $S$ is integral, then by the relative Picard scheme of $\pi$ we mean the relative Picard scheme of the restriction of $\pi$ to its smooth locus in $S$. A rational section of the relative Picard scheme will be equivalent to a regular section defined over the smooth locus of $\pi$.
\end{remark}

Following \citep[Remark 2.11]{kleiman}, we have an exact sequence \[ 0 \to \pic(S) \to \pic(\sD) \to \Pic(\sD/S)(S) \to \Br(S) \] where $\Pic(\sD/S)(S)$ is the group of sections of the relative Picard scheme. Since the Brauer group of any smooth variety is torsion, and the Brauer group of a curve is trivial by Tsen's theorem \citep[Proposition 13.6]{milne}, we have the following corollaries.

\begin{corollary} \label{brauertorsion}
Let $\tau \in \mathrm{NS}(\sD/S)$. Let $\sigma:S \to \Pic^\tau(\sD/S)$ be a section. There is a natural number $m$ such that $\sigma^{\otimes m}$ comes from a line bundle on $\sD$.
\end{corollary}

\begin{corollary} \label{tsen}
Let $\pi: \sD \to S$ be as above. Suppose now that $S$ has dimension one. Let $\tau \in \mathrm{NS}(\sD/S).$ Let $\sigma:S \to \Pic^\tau(\sD/S)$ be a section. Then there is a line bundle $L$ on $\sD$ which gives rise to $\sigma$.
\end{corollary}

\section{Counterexamples}\label{counterexamples}

In this section, we show that the hypothesis on the dimension of the locus of non-integral divisors in theorem \ref{mainresult} is necessary.

\begin{example}
Take the complete linear system of conics in $\bP^2$. The reducible conics form a divisor in this linear system. There is certainly a rational section of the relative Picard scheme which assigns to a smooth conic $C$ the line bundle $\cO_{\bP^1}(1)$. This section cannot arise from a line bundle on $\bP^2$, since the restriction of any line bundle on $\bP^2$ to a smooth conic has even degree.
\end{example}

Since conics have genus 0, they are somewhat exceptional. We will next give counterexamples of arbitrarily high genus.

\begin{example}
Let $\phi:S \to \bP^2$ be a double cover branched over a very general curve of degree $2d \geq 6$. Then $S$ is a surface of Picard number 1, generated by the pullback of $\cO_{\bP^2}(1)$. Let $C$ be the preimage of a conic in $\bP^2$. Then $C$ is a hyperelliptic curve of genus $2d-1$. The dimension of the linear system $|C|$ is five and every member of $|C|$ is a double cover of a conic, since \[ H^0(\cO_S(2)) \cong H^0(\phi_* \cO_S(2)) \cong H^0(\cO_{\bP^2}(2) \oplus \cO_{\bP^2}(2-d)). \]

The preimages of the reducible conics form a divisor of reducible curves in this linear system. There is a rational section of the degree two component of the relative Picard scheme which sends each curve to the pullback of $\cO_{\bP^1}(1)$ by the hyperelliptic map to $\bP^1$. On the other hand, it is easy to check that there is no line bundle on $S$ which has intersection number 2 with $C$, so this rational section cannot come from a line bundle on $S$.
\end{example}

\section{Proof of the Franchetta Conjecture for Linear Systems}

Let $X$ be a smooth projective variety, and $|D|$ a linear system on $X$. The universal divisor $\sD$ over $|D|$ maps to $X$, so we can pull back any line bundle $L$ on $X$ to the universal divisor. This gives a rational section $\sigma_L$ of the relative Picard scheme $\Pic(\sD/|D|)$, and its image is contained in some component $\Pic^\tau(\sD/|D|)$.

We will let $|D|^{s}$ be the complement of the discriminant locus in $|D|$, and $\sD^s$ its preimage in $\sD$. By Bertini's theorem, if $|D|$ is basepoint-free, then $|D|^{s}$ is nonempty. For the rest of this section, we will assume that the hypotheses of Theorem \ref{mainresult} hold, i.e., that $|D|$ is basepoint-free and the locus of non-integral divisors has codimension at least two.

We first prove a much weaker analog of Theorem \ref{mainresult}.

\begin{lemma}The relative Picard group $\pic(\sD/|D|)$ is generated by $\Pic(X)$.
\end{lemma}
\begin{proof}
Since $|D|$ is basepoint-free, the natural map $\sD \to X$ realizes $\sD$ as a projective bundle over $X$, so its Picard group is the direct sum of $\pic(X)$ and the tautological quotient line bundle $\cO(1)$. The latter is the pullback of $\cO_{|D|}(1)$.
\end{proof}

We now show that a rational section $\sigma$ of $\Pic(\sD/|D|)$ comes from a line bundle on $X$. The first step to proving Theorem \ref{mainresult} will be to show this \emph{pointwise}.
%We now recall the following. % Try to find a reference

%\begin{lem}
%Let $\tilde{X}$ be the blowup of $X$ at a smooth subvariety. Then $\Pic(\tilde{X}) \cong \Pic(X) \oplus G$, where $G$ is the free Abelian group generated by the connected components of the exceptional locus.
%\end{lem}

\begin{proposition} \label{keyProp}
For a general $[D'] \in |D|^s$, we have $\sigma([D'])=[L_{D'}|_{D'}]$ for some $L_{D'} \in \Pic(S)$.
\end{proposition}
\begin{proof}
%Consider a general pencil of divisors containing $D'$. By the hypothesis on the codimension of the non-integral locus, every divisor in the pencil is integral. By proposition \ref{section_extend}, we know that the rational section of the relative Picard scheme is define on the entire pencil.

Consider a general pencil in $|D|$. Let $\tilde{X}$ be the total space of the pencil, and $\bP^1$ the base. Then $\tilde{X}$ is the blowup of $X$ at the scheme-theoretic base locus of the pencil, which is smooth by Bertini's theorem. We have \[\pic(\tilde{X})\cong \pic(X) \oplus \bigoplus_i \Z E_i,\] where the $E_i$ are the connected components of the exceptional locus. Let $D_p$ denote the fiber of the natural map $\tilde{X} \to \bP^1$ over a point $p$. The class of $D_p$ in $\pic(\tilde{X})$ is $D-E$, where $E=\sum E_i$ is the exceptional divisor.

Let $C \subset \bP^1$ be an affine curve contained in the complement of the discriminant locus. By Proposition \ref{section_extend}, $\sigma$ can be uniquely extended to all of $C$. Let $D_C$ be the preimage of $C$ in $\tilde{X}$. All the fibers of the map $\tilde{X} \to \bP^1$ are integral by hypothesis, so \[\pic(D_C) \cong \pic(\tilde{X})/(D-E) \cong \pic(X) \oplus \Z E_i/(D-E)\] by the exact sequence of divisor class groups for an open subset \citep[II, 6.5]{hartshorne}.

By Corollary \ref{tsen}, $\sigma$ comes from a line bundle $\tilde{L}$ on $D_C$. By Corollary \ref{brauertorsion}, there is an integer $m$ such that $\sigma^{\otimes m}$ comes from some line bundle on $X$, from which it follows that $m\tilde{L}$ is in the image of $\pic(X)$. The following group-theoretic lemma will allow us to conclude that $\tilde{L}$ itself is in the image of $\pic(X)$.
\end{proof}

Let $A$ be an abelian group and let \[B = A \oplus \bigoplus_{i=1}^N E_i.\] Let $P$ be a partition of $N$ into $k$ parts $P_1, \dots, P_k$. For $1 \leq i \leq k$, let  \[r_i = a_i - \sum_{j \in P_i} E_j\] be elements of $B$. Let $C$ be the quotient group
$B/ \langle r_1, \dots, r_k\rangle$. Let $\tilde{A}$ denote the image of $A$ in $C$.

\begin{lemma}\label{groupLemma}
The subgroup $\tilde{A}$ is saturated in $C$, i.e., if there exists $c \in C$ with $nc \in \tilde{A}$ for some $n>0$, then $c \in \tilde{A}$. 
\end{lemma}
\begin{proof}
Let $\tilde{c}=a+\sum \beta_i E_i$ be a lift of $c$.  Let $R$ be the subgroup generated by the $r_i$. By assumption, we know that $na-a'+\sum_i \sum_{j \in P_i} n\beta_j E_j \in R$ for some $a' \in A$. For fixed $i$, the $n \beta_j$ must all be the same number $n \alpha_i$ for all $j \in P_i$, so the $\beta_j$ must all be equal to $\alpha_i$ for $j \in P_i$. We then get \[\tilde{c}+\sum_{i=1}^k \alpha_i r_i \in A\] which maps to $c \in C$.
\end{proof}

\begin{remark}\label{relaxintegrality}
The integrality hypothesis only comes to play in the lemma above, where it gives us that the relation is of the form $D-\sum E_i$ where each $E_i$ occurs precisely once. If for a general pencil, there is precisely one non-integral fiber, this fiber is reduced, and the base points of the pencil are contained in the smooth points of this fiber, then the hypotheses of this lemma will still be satisfied, and every other part of the argument will go through.
\end{remark}

We now want to deduce Theorem \ref{mainresult} from Proposition \ref{keyProp}. First, we will show that the $L_{D'}$ can all be chosen in the same connected component of $\Pic(X)$.

Let $\tau$ be the element of $\NS(\sD^s/|D|^s)$ corresponding to the component containing $\sigma$. We have a natural map $\NS(X) \to \NS(\sD^s/|D|^s)$. Let $T$ be the preimage of $\tau$ under this map. For each $\tau' \in T$, we have a restriction map \[\Pic^{\tau'}(X) \times |D|^s \to \Pic^\tau(\sD^s/|D|^s).\] Each of these maps is proper over $|D|^s$ since \[\Pic^{\tau'}(X) \times |D|^s \to \Pic^\tau(\sD^s/|D|^s)\] is proper over $|D|^s$ (since $\Pic^{\tau'}(X)$ is proper) and the map \[\Pic^\tau(\sD^s/|D|^s) \to |D|^s\] is proper. In particular, each of the restriction maps has a closed image. This means that the preimage of each $\Pic^{\tau'}$ in $|D|^s$ under $\sigma$ is a closed set. By Proposition \ref{keyProp}, the union of these closed sets is dense in $|D|^s$. By Severi's theorem of the base, $T$ is a countable set since $\NS(X)$ is. Therefore, one of these closed sets is all of $|D|^s$, or in other words, the image of $\sigma$ is contained in the image of $\Pic^{\tau'}$ for some fixed $\tau' \in \mathrm{NS}(X)$.
%, so by pulling back by $\sigma$  the images of the restriction maps for each $\tau' \in T$, we see that $|D|^s$ is a countable union of closed subvarieties, so one of them must be all of $|D|^s$, and hence there must be some $\tau'$ such that the image of $\sigma$ is contained in the image of $\Pic^{\tau'}(X)$.

Pick $\overline{L} \in \Pic^{\tau'}(X)$, and let $\sigma_{\overline{L}}$ be the corresponding section of $\Pic^\tau(\sD^s/|D|^s)$. By considering $\sigma-\sigma_{\overline{L}}$, we may as well assume that $\tau$ and $\tau'$ are both 0.

The restriction map \[r:\Pic^0(X) \times |D|^s \to \Pic^0(\sD^s/|D|^s)\] is a morphism of abelian schemes over $|D|^s$ which preserves 0. In particular, it is a group homomorphism, since this is true on each fiber \citep[II,Corollary 1]{mumford}. Let $K$ be the kernel of $r$. Let \[\pi:\Pic^0(X) \times |D|^s \to |D|^s\] be the projection onto the second factor. We would expect the kernel of this map to be constant over an open subset by semicontinuity. We now prove this.

\begin{lemma}
There is a nonempty open set $U \subset |D|^s$ such that $K \cap \pi^{-1}(x) \subset \Pic^0(X)$ is constant for $x \in U$.
\end{lemma}
\begin{proof}
By Grothendieck's theorem of generic flatness \citep[Theorem 14.4]{eisenbud}, there is a nonempty open set $V \subset |D|^s$ such that $K_x=K \cap \pi^{-1}(x) \subset \Pic^0(X)$ is a flat family of closed subvarieties of $\Pic^0(X)$. We will now restrict our attention to $V$.

Consider the component of the Hilbert scheme of closed subvarieties of $\Pic^0(X)$ which contains $K_x$. We want to show that no nearby point in the Hilbert scheme gives a subgroup of $\Pic^0(X)$. Basically, we show that all nearby points of the Hilbert scheme arise from translating each component of $K_x$ independently, and that if a nearby translate is not equal to $K_x$, then it cannot be a subgroup.

The tangent space to the corresponding point of the Hilbert scheme is given by $H^0(N_{K_x/\Pic^0(X)})$. Since $K_x$ is a closed subgroup, this normal bundle is a trivial bundle of rank equal to the codimension $c$ of $K_x$ in $\Pic^0(X)$. Therefore,  \[h^0(N_{K_x/\Pic^0(X)})=cn,\] where $n$ is the number of components of $K_x$.

We will now construct a flat family of embedded deformations of $K_x$ in $\Pic^0(X)$ such that its base dominates this component of the Hilbert scheme. Assume first that $K_x$ is connected. We note that we can identify the vector space $N_{K_x/\Pic^0(X),0}$ with $T_0(\Pic^0(X)/K_x)$. Let \[K'_x \subset \Pic^0(X) \times \Pic^0(X)\] be such that \[\pi_2^{-1}(\{a\})=K_x+a,\] i.e.~$K_x$ translated by $a$. This is just the universal family of translates of $K_x$. There is an induced map from $\Pic^0(X)$, considered as the base of this family, to the Hilbert scheme of subschemes of $\Pic^0(X)$, and the differential of this map at 0 is given by the natural map \[T_0 \Pic^0(X) \to T_0(\Pic^0(X)/K_x) \cong H^0(N_{K_x/\Pic^0(X)}),\] which is certainly surjective. Moreover, the kernel of this map consists of directions in which $a \in K_x$, or equivalently, directions in which $0 \in K_x+a$. In particular, any point near to $K_x$ but not equal to it, cannot be a subgroup, since it will not contain 0.

If $K_x$ is not connected, the above argument shows that any deformation preserving the subgroup structure must fix $K_x^0$, the connected component of the identity. This means we can identify deformations of $K_x$ with deformations of its image in $\Pic^0(X)/K_x^0$. But any such deformation preserving the subgroup structure must be contained in the torsion of $\Pic^0(X)/K_x^0$, which cannot happen for a nontrivial deformation.
\end{proof}

Let $K_0$ be $K_x$ for the $x$ in the $U$ of the above lemma. We get a birational factorization \[\Pic^0(X) \times |D|^s \to \Pic^0(X)/K_0 \times |D|^s \dashedrightarrow \Pic^0(\sD^s/|D|^s)\] where the last arrow is a rational map which is birational onto its image. We know that $\sigma$ is contained in the closure of the image of this last map, and $\sigma$ is defined for all points of $|D|^s$ by proposition \ref{section_extend}, so we see that $\sigma$ factors birationally to give a map \[|D|^s \dashedrightarrow \Pic^0(X)/K_0.\] Since $|D|^s$ is an open subvariety of projective space and $\Pic^0(X)/K_0$ is an abelian variety, this map must be constant. We can therefore find an element $L'$ of $\Pic^0(X)$ such that $\sigma$ and $\sigma_{L'}$ agree on a dense open subset of $|D|^s$, and hence agree everywhere.

We have now completed the proof of our first theorem.
\bigthmone*

Following the logic of Remark \ref{relaxintegrality}, we can actually prove something slightly stronger, which will be crucial for our proof of Theorem \ref{hitchinresult}.

\begin{corollary} \label{generalresult}
If $|D|$ is a basepoint-free linear system on $X$, and in a general pencil, there is precisely one non-integral fiber, this fiber is reduced, and the base locus of the pencil is disjoint from the singular locus of this fiber, then $\pic(X)$ surjects onto the Mordell-Weil group of the relative Picard scheme of the linear system.
\end{corollary}

\section{Higgs Bundles}

In this section, we will begin by reviewing the theory of Higgs bundles, and then use the machinery of the previous sections to determine all the sections of the Hitchin fibration. Except when stated otherwise, all material in this section can be found in \citep{nitsure} and \citep{simpson}. Through this section, $\Sigma$ will be a smooth projective curve of genus $g \geq 2$ with canonical bundle $K_\Sigma$.

\begin{definition}
A Higgs bundle on $\Sigma$ is a vector bundle $E$ on $\Sigma$ together with a map $\phi: E \to E \otimes K_\Sigma$, called the Higgs field. We define a Higgs bundle to be (semi)stable if for any nonzero proper subsheaf $F \subset E$ of strictly smaller rank which is preserved by $\phi$, $\mu(F)<\mu(E)$ (or $\mu(F) \leq \mu(E)$), where $\mu(F)$ is the usual Mumford slope $\deg(F)/\rk(F)$.
\end{definition}

\begin{theorem}
There is a moduli space of (S-equivalence classes of) semistable Higgs bundles of given rank $r$ and degree $d$, which we will denote $\cM(r,d)$.
\end{theorem}
\begin{proof}
This is \citep[Theorem 5.10]{nitsure}.
\end{proof}

Given any Higgs bundle of rank $r$, its characteristic polynomial is the element of the vector space \[B(r)=\bigoplus_{i=1}^r H^0(\Sigma,K_\Sigma^{\otimes i}) t^i\] defined by $\sum (-1)^i \tr(\phi^i) t^i$. Sending a semistable Higgs bundle to its characteristic polynomial gives a proper morphism \[H:\cM(r,d) \to B(r),\] called the \emph{Hitchin fibration} \citep[Theorem 6.1]{nitsure}. The vector space $B(r)$ is called the \emph{Hitchin base}. The general fibers of $H$ are abelian varieties, which we will later identify with the Jacobian of certain curves, called spectral curves. The rest of this paper is devoted to studying the Mordell-Weil group of this fibration.

To apply the machinery of the previous sections, we need to think of Higgs bundles as torsion sheaves on $T^*\Sigma$, the total space of the cotangent bundle of $\Sigma$. Note that the Higgs field $\phi:E \to E \otimes K_\Sigma$ gives $E$ the structure of a module over $\Sym(K_\Sigma^*)$, so it gives rise to a sheaf $\tilde{E}$ on \[\Spec(\Sym(K_\Sigma^*)) \cong T^*\Sigma\] which has pure one-dimensional support, called the \emph{spectral curve}.

Given any sheaf $\cF$ on $T^*\Sigma$ with pure one-dimensional support, we can push it forward to $\Sigma$ to get a vector bundle $F$ on $\Sigma$ with a map $F \to F \otimes K_\Sigma$. This correspondence gives rise to an isomorphism between the moduli space of Higgs bundles, $\cM(r,d)$, and the moduli space of pure one-dimensional semistable sheaves on $T^*\Sigma$ with determinant $r \pi^* K_\Sigma$ and Euler characteristic $d-r(g-1)$. The Hitchin fibration in this language is the map which sends a sheaf to its Fitting support, which is an element of the linear system $|r \pi^*K_\Sigma|$.

To compactify this moduli space, we can instead consider pure one-dimensional semistable sheaves on $X=\bP(\cO \oplus K_\Sigma)$, the projective completion of $T^*\Sigma$. In this case, the spectral curves will be elements of the linear system $\cO(r)$, which restricts to $\pi^*(K_\Sigma^r)$ on $T^*\Sigma$. This compactification adds new spectral curves which meet $D_\infty$, the divisor at infinity on $X$.

If we consider the locus of torsion sheaves with smooth, connected support, we get an open subscheme of $\cM(r,d)$. This subscheme is isomorphic to a component of the relative Picard scheme of the linear system $|\cO(r)|$, which is irreducible. This follows because if the Fitting support of a pure one-dimensional sheaf $\cF$ is a smooth connected curve $C$, then $\cF$ is the pushforward of a line bundle on $C$. Since this component of the Higgs moduli space is the only one which can dominate the Hitchin base, any rational section of the Hitchin fibration will be a rational section of the relative Picard scheme. We now have enough to begin the proof of Theorem \ref{hitchinresult}.

\bigthmtwo*

Given any line bundle on $\Sigma$, we clearly get a rational section of the relative Picard scheme, and hence an element of the Mordell-Weil group of the Hitchin fibration. We begin by describing these sections in terms of Higgs bundles.

First consider the case where the line bundle is $\cO_\Sigma$. Consider the vector bundle \[ E=\bigoplus_{i=1}^r K_\Sigma^{-i}. \] Given an element $\sigma$ of \[ B(r) \cong \bigoplus_{i=1}^r K_\Sigma^i,\] we get a natural map $K_\Sigma^{-r} \to E \otimes K_\Sigma$. Together with the identity maps $K_\Sigma^{-i} \to K_\Sigma^{-i}$ with $i$ between 1 and $r-1$, this gives us a map $E \to E \otimes K_\Sigma$, i.e., a Higgs field on $E$. The characteristic polynomial of this Higgs bundle will be $\sigma$. If we start with a different line bundle $L$, the Higgs bundle we get will be $E \otimes L$ with Higgs field $\phi \otimes \id_L$.

%Hom(E, E \otimes K)=H^0(\bigoplus K_\Sigma ^{-i} \otimes \bigoplus K_\Sigma^{i+1}
% On T^* \Sigma functions are given by H^0(O + \pi^*K + K^2+ ...) We have an element \sigma of H^0(K+K^2+...+K^r) which we mod out by. We have a tautological section of \pi^*K

We want to show that all sections of the Hitchin fibration are of this form. Ideally, we could then apply Theorem \ref{mainresult} directly, but unfortunately, the locus of non-integral curves in the linear system $|\cO(r)|$ has a divisorial component consisting of those curves which meet $D_\infty$. We first show that if this is the only divisorial component, then the Mordell-Weil group will consist entirely of these sections.

\begin{definition}
We will say that $\Sigma$ is $r$-\emph{good} if this is the only divisorial component of the non-integral locus of $\cO(r)$.
\end{definition}

We will begin by showing that if $\Sigma$ is $r$-good, the non-integral curves behave reasonably well.

\begin{proposition}
If $Sigma$ is $r$-good, then in a general pencil of curves in the linear system $|\cO(r)|$, there is a unique non-integral fiber which is reduced and reducible. Furthermore, the base locus of this pencil is disjoint from the singular locus.
\end{proposition}

\begin{proof}

We begin by showing that the locus of curves in the linear system which meet $D_\infty$ is a hyperplane, and every such curve must contain $D_\infty$ as a component. Pick any point $p \in D_\infty$. The locus of curves which meet $p$ is clearly a hyperplane. But since the pullback of $\cO(r)$ to $D_\infty$ is trivial, any curve which meets $p$ must contain $D_\infty$ as a component. It follows that in a general pencil, there is a unique non-integral curve, namely the curve corresponding the intersection of the line of the pencil and the hyperplane in $|\cO(r)|$ of curves containing $p$.

Since $\cO(r)$ is not a multiple of the divisor at infinity, which is irreducible, any curve in this linear system which contains $D_\infty$ must be reducible. We now show that the general such curve is reduced, for which purpose it suffices to exhibit a single such curve.

To do this, take the union of $D_\infty$, the pullback of a general section of $L$ on $\Sigma$, and a general section of $\cO(r-1)$. This is clearly a reduced curve containing $D_\infty$, and it is easy to verify that this reducible curve lies in the linear system $|\cO(r)|$.

For the last remark, it again suffices to give a single pencil with this property. Take the curve of the previous paragraph. Since the singular locus of this curve is finite, a general member of the basepoint-free linear system $|\cO(r)|$ misses it, and so the pencil these curves span has the desired property.
\end{proof}

\begin{proposition} \label{hitchingoodhkr}
If $\Sigma$ is $r$-good, then the Mordell-Weil group of the Hitchin fibration is generated by line bundles on $\Sigma$.
\end{proposition}
\begin{proof}
This follows from Corollary \ref{generalresult} together with the previous proposition.
\end{proof}

Now that we have a good idea of what happens when $\Sigma$ is $r$-good, we want to determine when this is the case. We begin by giving sufficient conditions in terms of the dimensions of certain cohomology groups.

\begin{proposition}\label{cohogoodness}
The curve $\Sigma$ is $r$-good if for $r_1+r_2=r$ with $r_i$ positive integers such that $r_1 \leq r_2$, we have \[ \sum_{i=1}^{r_1} h^0(K_\Sigma^i) <\sum_{i=r_2+1}^r h^0(K_\Sigma^i)-1. \]
\end{proposition}
\begin{proof}
It is clear from the definition of goodness and the discussion in the proof of the Proposition \ref{hitchingoodhkr} that we need to show that in the space of sections of $\cO(r)$ which do not meet $D_\infty$, the locus of non-integral sections has codimension strictly bigger than one. Any section not vanishing along $D_\infty$ comes from taking the product of elements of $H^0(\cO(r_1))$ and $H^0(\cO(r_2))$ which do not vanish along $D_\infty$. We can clearly assume here that $r_1 \leq r_2$. %with $r_1 \leq r_2$, $\dim |\cO(r_1)| +\dim |\cO(r_2)| < \dim |\cO(r)|-1$.

We want to show that $\dim |\cO(r_1)| +\dim |\cO(r_2)| < \dim |\cO(r)|-1$, since this will bound the dimension of the locus of non-integral curves. We can calculate the dimensions of the three linear systems, so this inequality becomes \[ \sum_{i=1}^{r_1} h^0(K_\Sigma^i) + \sum_{i=1}^{r_2} h^0(K_\Sigma^i) <\sum_{i=1}^r h^0(K_\Sigma^i)-1. \] This inequality is clearly equivalent to the one given in the statement of the theorem.
\end{proof}

%\begin{corollary}
%For any $L$, $(L,1)$ is good.
%\end{corollary}

\begin{corollary}
The curve $\Sigma$ is always $r$-good unless $r=g=2$.
\end{corollary}
\begin{proof}
If $r=2$ and $g>2$, then we get $r_1=r_2=1$. In this case, the conclusion follows immediately from Riemann-Roch. If $r>2$, we can assume that $r_2>1$. It then suffices to show that $h^0(K_\Sigma^i)<h^0(K_\Sigma^{i+r_2})-1$ for each $i$ between $1$ and $r_1$, and this follows from Riemann-Roch.
\end{proof}

We have now proved that $\Sigma$ is $r$-good unless $r=g=2$, but Proposition \ref{hitchingoodhkr} then implies that the Mordell-Weil group of the Hitchin fibration consists only of sections coming from line bundles on $\Sigma$. This completes the proof of our second theorem.

\subsection{Twisted Higgs Bundles}

We can extend this result by considering $L$-twisted Higgs bundles, where $L$ is a line bundle on $\Sigma$.

\begin{definition}
An $L$-twisted Higgs bundle is a vector bundle $E$ with a map $\phi:E \to E \otimes L$.
\end{definition}

The basic theory of Higgs bundles works essentially unchanged for twisted Higgs bundles if $L$ is basepoint-free of positive degree. We have a moduli space $\cM_L(r,d)$ of semistable $L$-twisted Higgs bundles with a Hitchin fibration. Let $X=\bP(\cO \oplus L)$. The space $\cM_L(r,d)$ is still sandwiched between a component of the relative Picard scheme of $|\cO_X(r)|$ and the moduli space of pure one-dimensional sheaves on $X$. We also want to extend our definition of $r$-goodness to the twisted case.

\begin{definition}
The line bundle $L$ is $r$-good if the only divisorial component of the locus of non-integral curves in $\cO_X(r)$ consists of curves meeting $D_\infty$.
\end{definition}

The proof of Proposition \ref{cohogoodness} gives us the following.

\begin{proposition}
The line bundle $L$ is $r$-good if for $r_1+r_2=r$ with $r_i$ positive integers such that $r_1 \leq r_2$, we have \[ \sum_{i=1}^{r_1} h^0(L^i) <\sum_{i=r_2+1}^r h^0(L^i)-1. \]
\end{proposition}

\begin{corollary}
If $g \geq 2$, $d=\deg(L) \geq 2g-2$, and $L \neq K_\Sigma$, then $L$ is $r$-good for any $r>1$.
\end{corollary}
\begin{proof}
Under these hypotheses, we have $h^0(L^i)=\chi(L^i)$. By Riemann-Roch, we have \[h^0(L^{r_2+i})-h^0(L^i)=d(r_2+i)-g+1-di+g-1=dr_2.\] We know that $r_2 \geq 1$ and $d \geq 2$, so \[\sum_{i=r_2+1}^r h^0(L^i)-1-\sum_{i=1}^{r_1} h^0(L^i) \geq 2r_1 \geq 2. \]
\end{proof}

We can now prove a version of Theorem \ref{hitchinresult} for twisted Higgs bundles.

\begin{theorem}
Let $L$ be a line bundle on $\Sigma$ which is $r$-good, e.g., let $g \geq 2$ and $\deg(L) \geq 2g-2$ with $L \neq K_\Sigma$. Then all rational sections of the Hitchin fibration for $\cM_L(r,d)$ come from line bundles on $L$.
\end{theorem}
\begin{proof}
Given that $L$ is $r$-good, the proof of Theorem \ref{hitchinresult} works unchanged in this setting.
\end{proof}

\bibliographystyle{plain}
\bibliography{divisor-franchetta-bibliography}

\end{document}